\documentclass[11pt]{article}
\usepackage{epic,latexsym,amssymb}
\usepackage{color}
\usepackage{tikz}
\usepackage{amsfonts,epsf,amsmath}

\newtheorem{thm}{Theorem}
\newtheorem{lem}{Lemma}

\newcommand{\qed}{$\Box$}

\newcommand{\proof}{\noindent\textbf{Proof. }}

\let\oldenumerate\enumerate
\renewcommand{\enumerate}{
  \oldenumerate
  \setlength{\itemsep}{0pt}
  \setlength{\parskip}{0pt}
  \setlength{\parsep}{0pt}
}

\begin{document}

\title{Independence and Connectivity of Connected Domination Critical Graphs}
%\subtitle{ The Arahato Samma Sumbudh Thasa}

\author{P. Kaemawichanurat$^1$\thanks{This work was funded by Development and Promotion of Science Technology Talents (DPST) Research Grant 031/2559.} \, and \, L. Caccetta$^2$
\\ \\
$^1$Theoretical and Computational Science Center \\
and Department of Mathematics, Faculty of Science,\\
King Mongkut's University of Technology Thonburi \\
Bangkok, Thailand\\
$^2$Western Australian Centre of Excellence in \\
Industrial Optimisation(WACEIO),\\
Department of Mathematics and Statistics, \\
Curtin University, Perth, WA 6845, Australia\\
\small \tt Email: $^{1}$pawaton.kae@kmutt.ac.th, $^{2}$L.Caccetta@exchange.curtin.edu.au}

\date{}
\maketitle

\begin{abstract}
A graph $G$ is said to be $k$-$\gamma_{c}$-critical if the connected domination number $\gamma_{c}(G) = k$ and $\gamma_{c}(G + uv) < k$ for every $uv \in E(\overline{G})$. Let $\delta, \kappa$ and $\alpha$ be respectively the minimum degree, the connectivity and the independence number. In this paper, we show that a $3$-$\gamma_{c}$-critical graph $G$ satisfies $\alpha \leq \kappa + 2$. Moreover, if $\kappa \geq 3$, then $\alpha = \kappa + p$ if and only if $\alpha = \delta + p$ for all $p \in \{1, 2\}$. We show that the condition $\kappa + 1 \leq \alpha \leq \kappa + 2$ is best possible to prove that $\kappa = \delta$. By these result, we conclude our paper with an open problem on Hamiltonian connected of $3$-$\gamma_{c}$-critical graphs.
\end{abstract}


%\vskip 1 cm
%\newpage
\section{\bf Introduction}
Our basic graph theoretic notation and terminology follows that of Bondy and Murty\cite{BM}. Thus $G$ denotes a finite graph with vertex set $V(G)$ and edge set $E(G)$. For $S \subseteq V(G)$, $G[S]$ denotes the subgraph of $G$ \emph{induced} by $S$. Throughout this paper all graphs are simple and connected. The \emph{open neighborhood} $N_{G}(v)$ of a vertex $v$ in $G$ is $\{u \in V(G) : uv \in E(G)\}$. Further, the \emph{closed neighborhood} $N_{G}[v]$ of a vertex $v$ in $G$ is $N_{G}(v) \cup \{v\}$. The \emph{degree} $deg_{G}(v)$ of a vertex $v$ in $G$ is $|N_{G}(v)|$. The \emph{minimum degree} of a graph $G$ is denoted by $\delta(G)$. $N_{S}(v)$ denotes $N_{G}(v) \cap S$ where $S$ is a vertex subset of $G$. A \emph{tree} is a connected graph that contains no cycle. A \emph{star} $K_{1, n}$ is a tree of $n + 1$ vertices which has $n$ vertices of degree $1$ and exactly one vertex of degree $n$. An \emph{independent set} is a set of pairwise non-adjacent vertices. The \emph{independence number} $\alpha(G)$ is the maximum cardinality of an independent set. For a connected graph $G$, a \emph{cut set} is a vertex subset $S \subseteq V(G)$ such that $G - S$ is connected. The minimum cardinality of a vertex cut set of a graph $G$ is called the \emph{connectivity} and is denoted by $\kappa(G)$. If $G$ has $S = \{a\}$ as a minimum cut set, then $a$ is a \emph{cut vertex} of $G$ and $\kappa(G) = 1$. A graph $G$ is \emph{$s$-connected} if $\kappa(G) \leq s$. When ambiguity occur, we abbreviate $\delta(G), \alpha(G)$ and $\kappa(G)$ to $\delta, \alpha$ and $\kappa$, respectively. A \emph{hamiltonian path} is a path containing every vertex of a graph. A graph $G$ is \emph{hamiltonian connected} if every two vertices of $G$ are joined by a hamiltonian path. Obviously, a graph of connectivity one is not hamiltonian connected. Moreover, a graph of connectivity two does not contain a hamiltonian path joining the two vertices in a minimum cut set. Thus we always focus on $3$-connected graphs when we study on the hamiltonian connected property of graphs. For a graph $G$ having a vertex subset $S$, the graph $G - S$ is obtained from $G$ by removing all vertices in $S$ and removing all edges that are incident to a vertex in $S$.
\vskip 5 pt

\indent For subsets $D, X \subseteq V(G)$, $D$ \emph{dominates} $X$ if every vertex in $X$ is either in $D$ or adjacent to a vertex in $D$. If $D$ dominates $X$, then we write $D \succ X$, further, we write $a \succ X$ when $D = \{a\}$. If $X = V(G)$, then $D$ is a \emph{dominating set} of $G$ and we write $D \succ G$ instead of $D \succ V(G)$. A smallest dominating set is call a $\gamma$\emph{-set}. The cardinality of a $\gamma$-set of $G$ is called the \emph{domination number} of $G$ and is denoted by $\gamma(G)$. A \emph{connected dominating set} of a graph $G$ is a dominating set $D$ of $G$ such that $G[D]$ is connected. If $D$ is a connected dominating set of $G$, we then write $D \succ_{c} G$. A smallest connected dominating set is call a $\gamma_{c}$\emph{-set}. The cardinality of a $\gamma_{c}$-set of $G$ is called the \emph{connected domination number} of $G$ and is denoted by $\gamma_{c}(G)$. A \emph{total dominating set} of a graph $G$ is a subset $D$ of $V(G)$ such that every vertex in $G$ is adjacent to a vertex in $D$. The minimum cardinality of a total dominating set of $G$ is called the \emph{total domination number} of $G$ and is denoted by $\gamma_{t}(G)$.
\vskip 5 pt

\indent A graph $G$ is $k$-$\gamma$-\emph{critical} if $\gamma(G) = k$ and $\gamma(G + uv) < k$ for each pair of non-adjacent vertices $u$ and $v$ of $G$. A $k$-$\gamma_{c}$-\emph{critical graph} and a $k$-$\gamma_{t}$-\emph{critical graph} are similarly defined.
\vskip 12 pt

\indent This paper focuses on the relationship of the connectivity and the independence number of $3$-$\gamma_{c}$-critical graphs. For related results on $3$-$\gamma$-critical graphs, Zhang and Tian\cite{ZT} showed that the independence number of every $3$-$\gamma$-critical graph does not exceed $\kappa + 2$ and, moreover if $\alpha = \kappa + 2$, then $\kappa = \delta$. By this result, the natural question to ask, on $3$-$\gamma_{c}$-critical graphs, is :
\vskip 8 pt

\indent What is the upper bound of the independence number of $3$-$\gamma_{c}$-critical graphs in term of the connectivity ?
\vskip 8 pt

\indent In this paper, we show that every $3$-$\gamma_{c}$-critical graph satisfies $\alpha \leq \kappa + 2$. Moreover, if $G$ is a $3$-$\gamma_{c}$-critical graph such that $\kappa + 1 \leq \alpha \leq \kappa + 2$, then $\kappa = \delta$ with only one exception. We show that the condition $\kappa + 1 \leq \alpha \leq \kappa + 2$ is best possible to establish that $\kappa = \delta$. Finally, Hamiltonian connected of $3$-$\gamma_{c}$-critical graphs with relationship between $\kappa$ and $\alpha$ are discussed in Section 4.
\vskip 5 pt

\section{\bf Preliminaries}
\noindent In this section, we state a number of results that we make use of in establishing our theorems. We begin with a result of Chv\'{a}tal and Erd\"{o}s\cite{CE} which gives hamiltonian properties of graphs according to the connectivity and the independence number.
\vskip 5 pt

\begin{thm}\cite{CE}\label{thm CE}
Let $G$ be an $s$-connected graph containing no independent set of $s$ vertices. Then $G$ is hamiltonian connected.
\end{thm}
\vskip 5 pt

\indent On $3$-$\gamma_{c}$-critical graphs, Ananchuen\cite{A} characterized such graphs with a cut vertex. For positive integers $n_{i}$ and $r \geq 2$, let $H = \cup^{r}_{i = 1}K_{1, n_{i}}$. For $1 \leq j \leq r$, let $c_{j}$ be the center of $K_{1, n_{j}}$ in $H$ and $w^{j}_{1}, w^{j}_{2}, ..., w^{j}_{n_{j}}$ the end vertices of $K_{1, n_{j}}$ in $H$. The graphs $G_{1}$ and $G_{2}$ are defined as follows. Set $V(G_{1}) = V(H) \cup \{x, y\}$ and $E(G_{1}) = E(\overline{H}) \cup \{xy\} \cup \{xw^{j}_{i}: 1 \leq i \leq n_{j}$ and $1 \leq j \leq r\}$. Set $V(G_{2}) = V(H) \cup \{x, y\} \cup U$ where $|U| \geq 1$ and $E(G_{2}) = E(\overline{H}) \cup \{xy\} \cup \{xw^{j}_{i}: 1 \leq i \leq n_{j}$ and $1 \leq j \leq r\} \cup \{uz: u \in U$ and $z \in V(G_{2}) \setminus \{x, y, u\}\}$. Ananchuen\cite{A} proved that :
\vskip 5 pt

\begin{thm}\label{thm A}\cite{A}
$G$ is a $3$-$\gamma_{c}$-critical graph with a cut vertex if and only if $G \in \{G_{1}, G_{2}\}$.
\end{thm}
\vskip 5 pt

\indent These following two lemmas, observed by Chen et al.\cite{CSM}, give fundamental properties of $3$-$\gamma_{c}$-critical graphs.
\vskip 5 pt

\begin{lem}\label{lem 1}\cite{CSM}
Let $G$ be a $3$-$\gamma_{c}$-critical graph and, for a pair of non-adjacent vertices $u$ and $v$ of $G$, let $D_{uv}$ be a $\gamma_{c}$-set of $G + uv$. Then
\begin {enumerate}
\item [(1)] $|D_{uv}| = 2$,
\item [(2)] $D_{uv} \cap \{u, v\} \neq \emptyset$ and
\item [(3)] if $u \in D_{uv}$ and $v \notin D_{uv}$, then $N_{G}(v) \cap D_{uv} = \emptyset$.
\end {enumerate}
\end{lem}
\vskip 5 pt

\begin{lem}\label{lem 2}\cite{CSM}
Let $G$ be a $3$-$\gamma_{c}$-critical graph and $I$ an independent set with $|I| = p \geq 3$. Then the vertices in $I$ can be ordered as $a_{1}, a_{2}, ..., a_{p}$ and there exists a path $x_{1}, x_{2},..., x_{p - 1}$ in $V(G) \setminus I$ with $\{a_{i}, x_{i}\} \succ_{c} G + a_{i}a_{i + 1}$ for $1 \leq i \leq p - 1$.
\end{lem}
\vskip 5 pt

\indent In the study of $3$-$\gamma_{t}$-critical graphs, Simmons\cite{S} proved that such graphs satisfy $\alpha \leq \delta + 2$. It was pointed out in Ananchuen\cite{A} that a $3$-$\gamma_{c}$-critical graph is also $3$-$\gamma_{t}$-critical and vice versa. Then every $3$-$\gamma_{c}$-critical graph satisfies $\alpha \leq \delta + 2$. In \cite{KCA}, we established the result of $3$-$\gamma_{c}$-critical graphs when $\alpha = \delta + 2$.
\vskip 5 pt

\begin{thm}\label{thm K}(\cite{KCA}, \cite{S})
Let $G$ be a $3$-$\gamma_{c}$-critical graph with $\delta \geq 2$. Then $\alpha \leq \delta + 2$. Moreover if  $\alpha = \delta + 2$, then $G$ contains exactly one vertex $x$ of degree $\delta$ and $G[N[x]]$ is a clique.
\end{thm}
\vskip 5 pt

%%%%%%%%%%%%%%%%%%%%%%%%%%%%%%%%%%%%%%%%%%%%%%%%%%%%%%%%%%%%%%%%%%%%%%%%%%%%%%%%%%%%%%%%%%%%%%%%%%%%%%%%%%%%%%%%%%%%%%%%%%%%%%%%%%%
\section{Main Results}
In this section we show that every $3$-$\gamma_{c}$-critical graph satisfies $\alpha \leq \kappa + 2$. If a $3$-$\gamma_{c}$-critical graphs contains a cut vertex, then, by Theorem \ref{thm A}, $\alpha = \kappa + 2$. Hence, throughout this paper, we focus on $2$-connected $3$-$\gamma_{c}$-critical graphs. Let $S$ be a minimum cut set and $C_{i}$ be the component of $G - S$ for $i = 1, 2, ..., m$. We, further, let $H_{1} = \cup^{\lfloor\frac{m}{2}\rfloor}_{i = 1}V(C_{i})$ and $H_{2} = \cup^{m}_{i = \lfloor\frac{m}{2}\rfloor + 1}V(C_{i})$. Let $I$ be a maximum independent set of $G$, $I_{i} = I \cap H_{i}$ and $|I_{i}| = \alpha_{i}$ for $i \in \{1, 2\}$. Then $I = I_{1} \cup I_{2} \cup (S \cap I)$. Let $|I_{1} \cup I_{2}| = p$. If $p \geq 3$, then, by Lemma \ref{lem 2}, the vertices in $I_{1} \cup I_{2}$ can be ordered as $a_{1}, a_{2}, ..., a_{s}$ and there exists a path $x_{1}, x_{2}, ..., x_{p - 1}$ in $V(G) \setminus (I_{1} \cup I_{2})$ such that $\{a_{i}, x_{i}\} \succ_{c} G + a_{i}a_{i + 1}$ for $1 \leq i \leq p - 1$. The following lemma gives the lower bound of $|S \setminus I|$.
\vskip 5 pt

\begin{lem}\label{lem P0}
If $p \geq 3$ and $|H_{1}|, |H_{2}| \geq 2$, then $x_{i} \in S \setminus I$, in particular, $|S \setminus I| \geq p - 1$.
\end{lem}
\proof
Suppose that $a_{i}, a_{i + 1} \in H_{j}$ for some $j \in \{1, 2\}$. To dominate $H_{3 - j}$, $x_{i} \in S \cup H_{3 - j}$. By the connectedness of $(G + a_{i}a_{i + 1})[D_{a_{i}a_{i + 1}}]$, $x_{i} \in S \setminus I$. We then suppose that $a_{i} \in H_{j}$ and $a_{i + 1} \in H_{3 - j}$. Since $|H_{3 - j}| \geq 2$, to dominate $H_{3 - j} \setminus \{a_{i + 1}\}$, $x_{i} \in H_{3 - j} \cup S$. Similarly, $x_{i} \in S \setminus I$. Thus $x_{1}, x_{2}, ..., x_{p - 1} \in S \setminus I$ implying that $|S \setminus I| \geq p - 1$. This completes the proof.
\qed
\vskip 5 pt

\indent We now ready to establish the upper bound of the independence number in term of the connectivity.
\vskip 5 pt

\begin{thm}\label{thm PM0}
The independence number of $3$-$\gamma_{c}$-critical graph does not exceed $\kappa + 2$.
\end{thm}
\proof
Suppose that $\alpha \geq \kappa + 3$. Thus $\alpha_{1} + \alpha_{2} + |S \cap I| \geq |S| + 3$. Therefore,
\begin{align}\label{eq p1}
p = \alpha_{1} + \alpha_{2} \geq |S \setminus I| + 3.
\end{align}

\indent We first show that $|H_{i}| \geq 2$ for all $i \in \{1, 2\}$. Suppose that $H_{i} = \{h\}$. Since $S$ is a minimum cut set, $N_{G}(h) = S$. We have that $\kappa = |S| = deg_{G}(h) \geq \delta$. Therefore, $\kappa = \delta$. So $\alpha \geq \kappa + 3 = \delta + 3$ contradicting Theorem \ref{thm K}. Hence, $|H_{1}|, |H_{2}| \geq 2$. Since $p \geq 3$, by Lemma \ref{lem P0}, $|S \setminus I| \geq p - 1$. Equation \ref{eq p1} implies that $p \geq |S \setminus I| + 3 \geq (p - 1) + 3 = p + 2$ a contradiction. Hence $\alpha \leq \kappa + 2$ and this completes the proof.
\qed

\indent To show the sharpness of the bound in Theorem \ref{thm PM0}, we introduce the class $\mathcal{G}_{1}(b_{0}, b_{1},$ $b_{2}, ..., b_{s - 1})$ where a graph $G$ in this class is defined as follows. Let $\{a_{0}, a_{1}, a_{2}, ..., a_{s}\}, B_{0},$ $ B_{1}... , B_{s - 1}$ be non-empty disjoint vertex sets $|B_{i}| = b_{i} \geq 1$ for all $0 \leq i \leq s - 1$. Let $b$ be a vertex in $B_{0}$. A graph $G$ in the class $\mathcal{G}_{1}(b_{0}, b_{1},$ $b_{2}, ..., b_{s - 1})$ can be constructed from $\{a_{0}, a_{1}, a_{2}, ..., a_{s}\}, B_{0}, B_{1}... , B_{s - 1}$ by adding edges according the join operations :
\begin{itemize}
  \item $a_{0} \vee (\{x\} \cup (\cup^{s - 1}_{j = 0}B_{j}))$,
  \item for $1 \leq i \leq s - 1$, $a_{i} \vee (\cup^{s - 1}_{j = 0}B_{j} \setminus B_{i})$,
  \item $a_{s} \vee (\cup^{s - 1}_{j = 1}B_{j})$,
  \item $x \vee (V(G) \setminus \{x, b, a_{s}\})$,
  \item adding edges so that the vertices in $\cup^{s - 1}_{j = 1}B_{j}$ form a clique and
  \item adding edges so that the vertices in $B_{0}$ form a clique.
\end{itemize}

\noindent We observe that $N_{G}(a_{s}) = \cup^{s - 1}_{i = 1}B_{i}, N_{G}(a_{0}) = \{x\} \cup (\cup^{s - 1}_{i = 0}B_{i}), N_{G}(x) = V(G) \setminus \{b, x, a_{s}\}$ and, for all $1 \leq i \leq s - 1$, $N_{G}(a_{i}) = \{x\} \cup (\cup^{s - 1}_{j = 0}B_{j} \setminus B_{i})$. We note further that if $B_{0} = \{b\}$, then $x$ is not adjacent to any vertex in $B_{0}$. A graph $G$ is illustrated in Figure 1.
\vskip 20 pt

\setlength{\unitlength}{0.7cm}
\begin{picture}(0, 11)
%Big G1

%vertices
\put(10.5, 0){\circle*{0.2}}
\put(4.5, 6){\circle*{0.2}}
\put(7.5, 6){\circle*{0.2}}
\put(9.5, 6){\circle*{0.2}}
\put(13.5, 6){\circle*{0.2}}
\put(15.5, 6){\circle*{0.2}}
\put(11, 9){\circle*{0.2}}

%oval
\put(10.5, 3){\oval(8, 1)}
\put(10.5, 9){\oval(2, 1)}

%line
\put(7.5, 6){\line(3, -1){6.5}}
\put(13.5, 6){\line(-5, -2){6.2}}

\put(13.5, 6){\line(1, 0){2}}
%a_{s} up
\put(10.5, 0){\line(3, 2){3.8}}
\put(10.5, 0){\line(-3, 2){3.8}}

%x up
\put(15.5, 6){\line(-2, 1){4.5}}
\multiput(15.5, 6)(-0.75,0.5){6}{\line(-3,2){0.6}}
\put(15.5, 6){\line(-6, 5){4.1}}
%x down
\put(15.5, 6){\line(-2, -3){1.65}}
\put(15.5, 6){\line(-3, -1){7.5}}

%a_{0} up
\put(4.5, 6){\line(3, 2){5.1}}
\put(4.5, 6){\line(2, 1){5.15}}
%a_{0} down
\put(4.5, 6){\line(2, -3){2}}
\put(4.5, 6){\line(4, -1){9.4}}

%a_{2} up
\put(9.5, 6){\line(0, 1){2.55}}
\put(9.5, 6){\line(2, 3){1.65}}
%a_{2} down
\put(9.5, 6){\line(-1, -1){2.5}}
\put(9.5, 6){\line(2, -1){4.5}}

%partition
\put(8.5, 2.5){\line(0, 1){1}}
\put(10.5, 2.5){\line(0, 1){1}}
\put(12.5, 2.5){\line(0, 1){1}}

\qbezier(15.5, 6)(12.5, 7.5)(9.5, 6)
\qbezier(15.5, 6)(11.5, 8)(7.5, 6)
\qbezier(15.5, 6)(10.5, 8.5)(4.5, 6)

\multiput(7.4, 3.5)(0,0.7){4}{\line(0,1){0.4}}
\multiput(7.6, 3.5)(0,0.7){4}{\line(0,1){0.4}}

\multiput(9.4, 3.5)(0,0.7){4}{\line(0,1){0.4}}
\multiput(9.6, 3.5)(0,0.7){4}{\line(0,1){0.4}}

\multiput(13.4, 3.5)(0,0.7){4}{\line(0,1){0.4}}
\multiput(13.6, 3.5)(0,0.7){4}{\line(0,1){0.4}}

%%%%%article%%%
\put(11.2, 3){$...$}
\put(11.2, 5.3){$...$}

\put(3.7, 6){$a_{0}$}
\put(6.8, 6){$a_{1}$}
\put(8.8, 6){$a_{2}$}
\put(12.3, 6){$a_{s - 1}$}
\put(9.7, 0){$a_{s}$}
\put(15.7, 6){$x$}

\put(7.5, 2.7){$B_{1}$}
\put(9.5, 2.7){$B_{2}$}
\put(13, 2.7){$B_{s - 1}$}
\put(10, 9.6){$B_{0}$}
\put(5, -1.5){\footnotesize\textbf{Figure 1 :} A graph $G$ in the class $\mathcal{G}_{1}(b_{0}, b_{1},$ $b_{2}, ..., b_{s - 1})$.}
%\put(17, 6){\footnotesize $a_{0}$}

%%%%%%%END Big G1%%%%%%%
\end{picture}
\vskip 70 pt

\begin{lem}\label{lem g1}
Let $G \in \mathcal{G}_{1}(b_{0}, b_{1},$ $b_{2}, ..., b_{s - 1})$. Then $G$ is a $3$-$\gamma_{c}$-critical graph.
\end{lem}
\proof
Let $b' \in B_{1}$. So $\{a_{0}, b', x\} \succ_{c} G$. Therefore $\gamma_{c}(G) \leq 3$. Suppose there exist $u, v \in V(G)$ such that $\{u, v\} \succ_{c} G$. If $\{u, v\} \cap (\cup^{s - 1}_{i = 1}B_{i}) \neq \emptyset$, then $|\{u, v\} \cap (\cup^{s - 1}_{i = 1}B_{i}|) = 1$ and $|\{u, v\} \cap \{a_{0}, a_{1}, ..., a_{s - 1}\}| = 1$ to dominate $B_{0}$. Without loss of generality, let $u \in B_{i}$ for some $i \in \{1, 2, ..., s - 1\}$. By connected and to dominate $b$, $v = \{a_{0}, a_{1}, ..., a_{s - 1}\} \setminus \{a_{i}\}$. Thus $\{u, v\}$ does not dominate $a_{i}$. Hence $\{u, v\} \cap (\cup^{s - 1}_{i = 1}B_{i}) = \emptyset$. This together with $G[\{u, v\}]$ is connected imply that either $\{u, v\} \subseteq \{x, a_{0}, a_{1}, ..., a_{s - 1}\} \cup B_{0}$ or $\{u, v\} \subseteq \{a_{s}\}$. Clearly the second case is not possible. Thus $\{u, v\} \subseteq \{x, a_{0}, a_{1}, ..., a_{s - 1}\} \cup B_{0}$ but $\{u, v\}$ does not dominate $a_{s}$, a contradiction. Hence $\gamma_{c}(G) = 3$.
\vskip 5 pt

\indent We next establish the criticality of $G$. For any pair of non-adjacent vertices $u, v \in V(G)$. Consider $G + uv$. If $\{u, v\} = \{x, b\}$, then $D_{uv} = \{x, b'\}$ for some $b' \in \cup^{s - 1}_{i = 1}B_{i}$. If $\{u, v\} = \{a_{s}, b''\}$ where $b'' \in B_{0}$, then $D_{uv} = \{b'', a_{0}\}$. If $\{u, v\} = \{a_{s}, a_{i}\}$ for $0 \leq i \leq s - 1$, then $D_{uv} = \{a_{i}, x\}$. If $\{u, v\} = \{a_{s}, x\}$, then $D_{uv} = \{a_{0}, x\}$. If $\{u, v\} = \{a_{i}, a_{j}\}$ or $\{u, v\} = \{b_{j}, a_{j}\}$  where $b_{j} \in B_{j}$, then $D_{uv} = \{a_{i}, b_{j}\}$. Therefore $G$ is $3$-$\gamma_{c}$-critical graph and this completes the proof.
\qed
\vskip 5 pt

\indent We see that a graph $G \in \mathcal{G}_{1}(b_{0}, 1, 1, ..., 1)$ has $\cup^{s - 1}_{i = 1}B_{i}$ as a minimum cut set with $s - 1$ vertices. Moreover, $G$ has $\{a_{0}, a_{1}, ..., a_{s}\}$ the maximum independent set with $s + 1$ vertices. Thus $\alpha = \kappa + 2$ and the bound in Theorem \ref{thm PM0} is best possible.
\vskip 10 pt

\indent In the following, we show that if a $3$-$\gamma_{c}$-critical graph satisfies $\kappa + 1 \leq \alpha \leq \kappa + 2$, then $\kappa = \delta$ with only one exception. We give a construction of the class $\mathcal{G}_{2}(3, 3)$ which a graph $G$ is this class is $3$-$\gamma_{c}$-critical satisfying $\alpha = \kappa + 1$ and $\kappa < \delta$. Let $J_{1}, J_{2}$ be two disjoint vertex sets such that$|J_{1}| = j_{1} \geq 3, |J_{2}| = j_{2} \geq 3$ and let $u_{1}, u_{2} \in J_{2}$. Moreover, let $w, w_{1}, v_{1}$ three isolated vertices. A graph $G$ in the class $\mathcal{G}_{2}(3, 3)$ can be constructed from $J_{1}, J_{2}, w, w_{1}$ and $v_{1}$ by adding edges according the join operations :

\begin{itemize}
  \item $v_{1} \vee J_{1}$,
  \item $w_{1} \vee (J_{1} \cup (J_{2} \setminus \{u_{1}\}) \cup \{w\})$,
  \item $w \vee (J_{1} \cup (J_{2} \setminus \{u_{2}\}) \cup \{w_{1}\})$,
  \item adding edges so that the vertices in $J_{1}$ form a clique and
  \item adding edges so that the vertices in $J_{2}$ form a clique.
\end{itemize}

\noindent We observe that $N_{G}(w_{1}) = J_{1} \cup (J_{2} \setminus \{u_{1}\}) \cup \{w\}$ and $N_{G}(w) = J_{1} \cup (J_{2} \setminus \{u_{2}\}) \cup \{w_{1}\}$. $G[J_{1}]$ and $G[J_{2}]$ are cliques. It is easy to see that $G$ is a $3$-$\gamma_{c}$-critical graph. A graph $G$ has $\{w, w_{1}\}$ as a minimum cut set and $\{u_{1}, w_{1}, v_{1}\}$ as a maximum independent set. Thus $\alpha = \kappa + 1$. Since $j_{1}, j_{2} \geq 3$, $\delta \geq 3$. Therefore $\delta > \kappa$. Before we give the next theorem, let us introduce a lemma which is the situation that might occur often in the following proofs.
\vskip 5 pt

\begin{lem}\label{lem w}
Let $G$ be a $3$-$\gamma_{c}$-critical graph and $b_{1}, b_{2}$ and $b_{3}$ be three non-adjacent vertices of $G$. Then $D_{b_{1}b_{2}}$ is either $\{b_{1}, b\}$ or $\{b_{2}, b\}$ for some $b \in V(G)$.
\end{lem}
\proof
Consider $G + b_{1}b_{2}$. Lemma \ref{lem 1}(1) implies that $|D_{b_{1}b_{2}}| = 2$, moreover, Lemma \ref{lem 1}(2) implies that $D_{b_{1}b_{2}} \cap \{b_{1}, b\} \neq \emptyset$. To dominate $b_{3}$, $D_{b_{1}b_{2}} \neq \{b_{1}, b_{2}\}$. Therefore $D_{b_{1}b_{2}}$ is either $\{b_{1}, b\}$ or $\{b_{2}, b\}$ for some $b \in V(G)$ and this completes the proof.
\qed

\begin{thm}\label{thm PM1}
Let $G$ be a $3$-$\gamma_{c}$-critical graph with a minimum cut set $S$ and a maximum independent set $I$. If $\kappa + 1 \leq \alpha, \kappa < \delta$ and $|S \setminus I| \leq 1$, then $G \in \mathcal{G}_{2}(3, 3)$.
\end{thm}
\proof
In view of Theorem \ref{thm PM0}, $\kappa + 1 \leq \alpha \leq \kappa + 2$. Recall that $p = |I_{1} \cup I_{2}|$. Thus $|S| + 1 \leq p + |S \cap I| \leq |S| + 2$. Therefore,
\begin{align}\label{eq p2}
|S \setminus I| + 1 \leq p \leq |S \setminus I| + 2.
\end{align}

\indent Suppose that $|H_{i}| = 1$ for some $i = 1$ or $2$. Let $H_{i} = \{h\}$. By the minimality of $S$, $N_{G}(h) = S$ which implies that
\begin{center}
$\kappa = deg_{G}(h) \geq \delta > \kappa$,
\end{center}
a contradiction. Clearly $|H_{i}| \geq 2$ for all $1 \leq i \leq 2$. Since $|S| \geq 2$ and $|S \setminus I| \leq 1$, it follows that $S \cap I \neq \emptyset$. We distinguish $2$ cases according to the number of vertices in $S \setminus I$.
\vskip 5 pt

\noindent \textbf{Case 1 :} $|S \setminus I| = 0$.\\
\indent Therefore, $S \subseteq I$. By Equation \ref{eq p2}, suppose first that $p = |S \setminus I| + 2$. Thus $p = 2$. Let $\{v_{1}, v_{2}\} = I_{1} \cup I_{2}$. Consider $G + v_{1}v_{2}$. Lemma \ref{lem 1}(2) implies that $D_{v_{1}v_{2}} \cap \{v_{1}, v_{2}\} \neq \emptyset$. To dominate $S \cap I$, $D_{v_{1}v_{2}} \neq \{v_{1}, v_{2}\}$. Without loss of generality, let $D_{v_{1}v_{2}} = \{v_{1}, v\}$ for some $v \in V(G)$. Since $S \subseteq I$ and $(G + v_{1}v_{2})[D_{v_{1}v_{2}}]$ is connected, it follows that $v$ is in the same set $H_{i}$ of $v_{1}$ for some $i \in \{1, 2\}$. Because $|H_{3 - i}| \geq 2$, $D_{v_{1}v_{2}}$ does not dominate $H_{3 - i} \setminus \{v_{2}\}$ contradicting Lemma \ref{lem 1}(1). Hence, $p = |S \setminus I| + 1 = 1$.
\vskip 5 pt

\indent Without loss of generality, let $\{v_{1}\} = I_{1} = I_{1} \cup I_{2}$. Because $|S| \geq 2$ and $S \subseteq I$, it follows that there exist $w_{1}, w_{2} \in S \cap I$. Consider $G + w_{1}w_{2}$. Clearly $w_{1}$ and $w_{2}$ are not adjacent to $v_{1}$. Lemma \ref{lem w} then yields that $D_{w_{1}w_{2}}$ is $\{w_{1}, w\}$ or $\{w_{2}, w\}$ for some $w \in V(G)$. Without loss of generality, let $D_{w_{1}w_{2}} = \{w_{1}, w\}$. To dominate $v_{1}$, $w \in H_{1}$. Thus $w_{1} \succ H_{2}$. Lemma \ref{lem 1}(3) gives that $ww_{2} \notin E(G)$. We now have that $w_{2}$ is not adjacent to $v_{1}$ and $w$ in $H_{1}$. Consider $G + w_{2}v_{1}$. Similarly, $|D_{w_{2}v_{1}} \cap \{w_{2}, v_{1}\}| = 1$ to dominate $w_{1}$. Thus $D_{w_{2}v_{1}} = \{v_{1}, u\}$ or $D_{w_{2}v_{1}} = \{w_{2}, u\}$ for some $u \in V(G)$. If $D_{w_{2}v_{1}} = \{v_{1}, u\}$, then $u \in S$ to dominate $H_{2}$. But $S \subseteq I$. Hence $(G + w_{2}v_{1})[D_{w_{2}v_{1}}]$ is not connected, a contradiction. Therefore, $D_{w_{2}v_{1}} = \{w_{2}, u\}$. To dominate $w$, $u \in H_{1} \cup S$. As $w_{2} \in S \subseteq I$, we must have by the connectedness of $(G + v_{1}w_{2})[D_{v_{1}w_{2}}]$ that $u \in H_{1}$. Thus $w_{2} \succ H_{2}, uw_{2} \in E(G)$ and $uv_{1} \notin E(G)$. Since $w_{1}w_{2} \notin E(G)$, $uw_{1} \in E(G)$. Let $y \in H_{2}$. Because $w_{1} \succ H_{2}$ and $w_{2} \succ H_{2}$, it follows that $yw_{1}, yw_{2} \in E(G)$. Consider $G + uy$. Clearly $u$ and $y$ are not adjacent to $v_{1}$. Lemma \ref{lem w} then gives that $D_{uy} = \{u, x\}$ or $D_{uy} = \{y, x\}$ for some $x \in V(G)$. Since $|H_{i}| \geq 2$ for $1 \leq i \leq 2$, $x \in S$ to dominate $H_{2} \setminus \{y\}$ when $D_{uy} = \{u, x\}$ and $H_{1} \setminus \{u\}$ when $D_{uy} = \{y, x\}$. Since $S \subseteq I$, $xv_{1} \notin E(G)$. Thus $D_{uy}$ does not dominate $v_{1}$ contradicting Lemma \ref{lem 1}(1). Thus Case 1 cannot occur.
\vskip 5 pt

\noindent \textbf{Case 2 :} $|S \setminus I| = 1$.\\
\indent Since $|S| \geq 2$, $S \cap I \neq \emptyset$. By Equation \ref{eq p2}, suppose first that $p = |S \setminus I| + 2 = 3$. By the pigeon' hole principle, there exists $i \in \{1, 2\}$ such that $|I_{i}| \geq 2$. Let $v_{1}, v_{2} \in I_{i}$ and $\{v_{3}\} = I_{1} \cup I_{2} \setminus \{v_{1}, v_{2}\}$. Consider $G + v_{1}v_{2}$. Clearly $v_{1}$ and $v_{2}$ are not adjacent to $v_{3}$. Lemma \ref{lem w} implies that $D_{v_{1}v_{2}} = \{v_{1}, u\}$ or $D_{v_{1}v_{2}} = \{v_{2}, u\}$ for some $u \in V(G)$. Without loss of generality, let $D_{v_{1}v_{2}} = \{v_{1}, u\}$. To dominate $H_{3 - i}$, $u \in S \cup H_{3 - i}$. By the connectedness of $(G + v_{1}v_{2})[D_{v_{1}v_{2}}]$, $\{u\} = S \setminus I$. So $u \succ \{v_{1}, v_{3}\}$ and $uv_{2} \notin E(G)$. This together with $|S \setminus I| = 1$, we have that $v_{2}$ is not adjacent to any vertex in $S$. Consider $G + v_{1}v_{3}$. Similarly, $D_{v_{1}v_{3}} = \{v_{1}, u'\}$ or $D_{v_{1}v_{3}} = \{v_{2}, u'\}$ for some $u' \in V(G)$. To dominate $v_{2}$, $u' \in H_{i}$. By the connectedness of $(G + v_{1}v_{3})[D_{v_{1}v_{3}}]$, $D_{v_{1}v_{3}} \subseteq H_{i}$. Thus $D_{v_{1}v_{3}}$ does not dominate $H_{3 - i} \setminus \{v_{3}\}$ contradiction Lemma \ref{lem 1}(1). Therefore, $p = |S \setminus I| + 1 = 2$.
\vskip 5 pt

\noindent \textbf{Subcase 2.1 :} $\alpha_{1} = 0$ and $\alpha_{2} = 2$.\\
\indent Let $v_{1}, v_{2} \in I_{2}$. Consider $G + v_{1}v_{2}$. To dominate $H_{1}$, it follows from Lemma \ref{lem w} that $D_{v_{1}v_{2}} = \{v_{1}, u\}$ or $D_{v_{1}v_{2}} = \{v_{2}, u\}$ for some $u \in V(G)$. Without loss of generality, let $D_{v_{1}v_{2}} = \{v_{1}, u\}$. Thus $u \succ H_{1} \cup (S \cap I)$ and $uv_{2} \notin E(G)$. So $\{u\} = S \setminus I$. Thus $v_{2}$ is not adjacent to any vertex in $S$. Let $v \in H_{1}$. Consider $G + vv_{1}$. Similarly, $D_{vv_{1}} = \{v, u'\}$ or $D_{vv_{1}} = \{v_{1}, u'\}$. To dominate $v_{2}$, $u' \in H_{2}$. By the connectedness of $(G + vv_{1})[D_{vv_{1}}]$, $D_{vv_{1}} = \{v_{1}, u'\}$. So $D_{vv_{1}}$ does not dominate $H_{1} \setminus \{v\}$ contradicting Lemma \ref{lem 1}(1). Consequently, Subcase 2.1 cannot occur. Similarly, $\alpha_{2} = 0$ and $\alpha_{1} = 2$ cannot occur.
\vskip 5 pt

\noindent \textbf{Subcase 2.2 :}  $\alpha_{1} = 1$ and $\alpha_{2} = 1$.\\
\indent Let $u_{1} \in I_{1}$ and $v_{1} \in I_{2}$. Consider $G + u_{1}v_{1}$. To dominate $I$, this implies by Lemma \ref{lem w} that $D_{u_{1}v_{1}} = \{u_{1}, w\}$ or $D_{u_{1}v_{1}} = \{v_{1}, w\}$ for some $w \in V(G)$. Without loss of generality, let $D_{u_{1}v_{1}} = \{u_{1}, w\}$. Thus $u_{1}w \in E(G)$ and $w \succ (S \cap I) \cup (H_{2} \setminus \{v_{1}\})$. Therefore $w \in S$. As $wu_{1} \in E(G)$ and $u_{1} \in I$, we must have $\{w\} = S \setminus I$.
\vskip 5 pt

\noindent \textbf{Claim 1 :} $|S \cap I| = 1$, in particular, $|S| = 2$.\\
\indent Suppose there exist $w_{1}, w_{2} \in S \cap I$. Consider $G + w_{1}w_{2}$. We see that $w_{1}$ and $w_{2}$ are not adjacent to $u_{1}$. This implies by Lemma \ref{lem w} that $D_{w_{1}w_{2}} = \{w_{1}, w'\}$ or $D_{w_{1}w_{2}} = \{w_{2}, w'\}$ for some $w' \in H_{1} \cup \{w\}$. In both cases, $D_{w_{1}w_{2}}$ does not dominate $\{v_{1}, u_{1}\}$ contradicting Lemma \ref{lem 1}(1). Therefore $|S \cup I| = 1$ and we settle Claim 1.
\vskip 5 pt

\indent Let $\{w_{1}\} = S \cap I$. As $w \succ S \cap I$, we must have $ww_{1} \in E(G)$. We have that $w_{1}$ is adjacent to neither $u_{1}$ in $H_{1}$ nor $v_{1}$ in $H_{2}$.
\vskip 5 pt

\noindent \textbf{Claim 2 :} $w_{1} \succ H_{1} \setminus \{u_{1}\}$ and $w_{1} \succ H_{2} \setminus \{v_{1}\}$.\\
\indent Consider $G + w_{1}u'$ where $u' \in \{u_{1}, v_{1}\}$. Suppose that $u' \in H_{i}$ for some $i \in \{1, 2\}$. We see that $w_{1}$ and $u'$ are not adjacent to the vertex in $\{u_{1}, v_{1}\} \setminus \{u'\}$. Lemma \ref{lem w} thus implies $D_{w_{1}u'} = \{w_{1}, v'\}$ or $D_{w_{1}u'} = \{u', v'\}$ for some $v' \in V(G) \setminus \{w_{1}, u'\}$. Suppose that $D_{w_{1}u'} = \{u', v'\}$. To dominate $H_{3 - i} \setminus \{u_{1}, v_{1}\}$ and by the connectedness of $(G + w_{1}u')[D_{w_{1}u'}]$, $v' \in S \setminus I$. So $v' = w$. But $ww_{1} \in E(G)$ contradicting Lemma \ref{lem 1}(3). Therefore, $D_{w_{1}u'} = \{w_{1}, v'\}$. To dominate $H_{3 - i}$, $v' \in H_{3 - i}$. So $w_{1} \succ H_{i} \setminus \{u'\}$ and thus establishing Claim 2.
\vskip 5 pt

\noindent \textbf{Claim 3 :} $w \succ H_{2} \setminus \{v_{1}\}$ and $w \succ H_{1} \setminus \{u_{2}\}$ for some $u_{2} \in H_{1} \setminus \{u_{1}\}$.\\
\indent Recall that $u_{1}w, w_{1}w \in E(G)$ and $w \succ H_{2} \setminus \{v_{1}\}$. By the assumption that $\kappa < \delta$, $v_{1}$ is adjacent to some vertex $v'$ in $H_{2}$. If $w \succ H_{1}$, then $\{w, v'\} \succ_{c} G$ contradicting $\gamma_{c}(G) = 3$. Thus there exists $u_{2} \in H_{1} \setminus \{u_{1}\}$ such that $wu_{2} \notin E(G)$. Consider $G + wu_{2}$. By the same arguments as Claim 2, $D_{wu_{2}} = \{w, w'\}$ where $w' \in H_{2} \setminus \{v_{1}\}$. Thus $w \succ H_{1} \setminus \{u_{2}\}$ and establishing Claim 3.
\vskip 5 pt

\noindent \textbf{Claim 4 :} $G[H_{i}]$ is a clique for all $1 \leq i \leq 2$.\\
\indent Suppose there exist $h_{1}, h_{2} \in H_{i}$ for some $i \in \{1, 2\}$ such that $h_{1}h_{2} \notin E(G)$. Consider $G + h_{1}h_{2}$. By Lemma \ref{lem 1}(1), $|D_{h_{1}h_{2}}| = 2$ and by Lemma \ref{lem 1}(2), $D_{h_{1}h_{2}}\cap \{h_{1}, h_{2}\} \neq \emptyset$. To dominate $H_{3 - i}$, $|D_{h_{1}h_{2}} \cap \{h_{1}, h_{2}\}| = 1$. As $(G + h_{1}h_{2})[D_{h_{1}h_{2}}]$ is connected, we must have that $|D_{h_{1}h_{2}} \cap S| = 1$. Thus $D_{h_{1}h_{2}} \cap S$ is either $\{w_{1}\}$ or $\{w\}$. If $D_{h_{1}h_{2}} \cap S = \{w_{1}\}$, then $D_{h_{1}h_{2}}$ does not dominate $\{u_{1}, v_{1}\}$, a contradiction. Hence $D_{h_{1}h_{2}} \cap S = \{w\}$. But this yields that $D_{h_{1}h_{2}}$ does not dominate $\{v_{1}, u_{2}\}$, a contradiction. So $G[H_{i}]$ is a clique for all $1 \leq i \leq 2$ and thus establishing Claim 4.
\vskip 5 pt

\indent In view of Claims 1, 2, 3 and 4, We have that $G \in \mathcal{G}_{2}(3, 3)$ and this completes the proof.
\qed
\vskip 5 pt

\begin{lem}\label{lem PM2}
Let $G$ be a $3$-$\gamma_{c}$-critical graph with a minimum cut set $S$ and a maximum independent set $I$ and $H_{i}$ be defined the same as in Theorem \ref{thm PM1}. If $|S \setminus I| \geq 2$ and $\kappa + 1 \leq \alpha$, then at least one of $H_{i}$ must be singleton.
\end{lem}
\proof
We may assume to the contrary that $|H_{i}| \geq 2$ for all $i = 1$ and $2$. As $p + |S \cap I| = \alpha \geq \kappa + 1 = |S \cap I| + |S \setminus I| + 1$, we must have
\begin{center}
$p \geq |S \setminus I| + 1 \geq 3$.
\end{center}

\noindent This implies by Lemma \ref{lem 2} that the vertices in $I_{1} \cup I_{2}$ can be order $a_{1}, a_{2}, ..., a_{p}$ and there exits a path $x_{1}, .., x_{p - 1}$ such that $\{a_{i}, x_{i}\} \succ_{c} G + a_{i}a_{i + 1}$ for all $1 \leq i \leq p - 1$. Lemma \ref{lem P0} yields that $x_{i} \in S \setminus I$ for all $1 \leq i \leq p - 1$. Thus
\begin{align}\label{eq SI}
p - 1 \leq |S \setminus I|.
\end{align}

\noindent This implies that
\begin{center}
$p - 1 = |S \setminus I|$.
\end{center}

\noindent That is $S \setminus I = \{x_{1}, .., x_{p - 1}\}$. Since $a_{i} \in I$, $x_{i} \succ S \cap I$ for all $1 \leq i \leq p - 1$.
\vskip 5 pt

\noindent \textbf{Claim 1 :} $G[\{x_{1}, x_{2}, ..., x_{p - 1}\}]$ is a clique.\\
\indent Consider $G + a_{i}a_{j}$ where $2 \leq i, j \leq p$. By Lemma \ref{lem 1}(2) and to dominate $I$, $|D_{a_{i}a_{j}} \cap \{a_{i}, a_{j}\}| = 1$. Without loss of generality, let $D_{a_{i}a_{j}} = \{a_{i}, v\}$ for some $v \in V(G)$. We show that $v \in \{x_{1}, x_{2}, ..., x_{p - 1}\}$. We first consider the case when $a_{i}, a_{j} \in H_{l}$ for some $l \in \{1, 2\}$. To dominate $H_{3 - l}$, $v \in H_{3 - l} \cup S$. By the connectedness of $(G + a_{i}a_{j})[D_{a_{i}a_{j}}]$, $v \in S \setminus I$. Thus $v \in \{x_{1}, x_{2}, ..., x_{p - 1}\}$. We now consider the case when $a_{i} \in H_{l}$ and $a_{j} \in H_{3 - l}$. Since $H_{3 - l} \geq 2$, $H_{3 - l} \setminus \{a_{j}\} \neq \emptyset$. To dominate $H_{3 - l}$ and by the connectedness of $(G + a_{i}a_{j})[D_{a_{i}a_{j}}]$, $v \in S \setminus I$. Thus $v \in \{x_{1}, .., x_{p - 1}\}$. In both cases, by Lemmas \ref{lem 1}(3) and \ref{lem 2}, we have that $v = x_{j - 1}$. But $a_{i}x_{i - 1} \notin E(G)$. It follows that $x_{i - 1}x_{j - 1} \in E(G)$ for all $2 \leq i, j \leq p$. Hence $G[\{x_{1}, x_{2}, ..., x_{p - 1}\}]$ is a clique and establishing Claim 1.
\vskip 5 pt

\indent Without loss of generality, let $\alpha_{1} \leq \alpha_{2}$.
\vskip 5 pt

\noindent \textbf{Case 1 :} $\alpha_{1} = 0$.\\
\indent Thus $a_{1}, a_{2}, ..., a_{p} \in H_{2}$ and $x_{i} \succ H_{1}$ for all $1 \leq i \leq p - 1$. Let $h \in H_{1}$. Therefore $h \succ S \setminus I$. Lemmas \ref{lem 2} and \ref{lem P0} imply that $a_{1} \succ S \setminus I$. Consider $G + ha_{1}$. Clearly $ha_{2}, a_{1}a_{2} \notin E(G)$. It follows from Lemma \ref{lem w} that $D_{ha_{1}} = \{h, w\}$ or $D_{ha_{1}} = \{a_{1}, w\}$ for some $w \in V(G)$. If $D_{ha_{1}} = \{h, w\}$, then $w \in H_{2} \cup S$ to dominate $H_{2}$. By the connectedness of $(G + ha_{1})[D_{ha_{1}}]$, $w \in S$. To dominate $a_{2}$, $w \in S \setminus I$. Thus $wa_{1} \in E(G)$ contradicting Lemma \ref{lem 1}(3). Therefore, $D_{ha_{1}} = \{a_{1}, w\}$. To dominate $H_{1} \setminus \{h\}$, $w \in H_{1} \cup S$. By the connectedness of $(G + ha_{1})[D_{ha_{1}}]$, $w \in S \setminus I$. Clearly $wh \in E(G)$ contradicting Lemma \ref{lem 1}(3). Thus Case 1 cannot occur.
\vskip 5 pt

\noindent \textbf{Case 2 :} $\alpha_{1} \geq 1$.\\
\indent Suppose that $\alpha_{1} = 1$ and $\{a_{p}\} = I_{1}$. Since $|H_{1}| \geq 2$, there exists $h' \in H_{1} \setminus \{a_{p}\}$. By Lemma \ref{lem 2}, and $a_{1}, a_{2}, ..., a_{p - 1} \in H_{2}$, $x_{i}h' \in E(G)$ for $1 \leq i \leq p - 1$. Therefore, $h' \succ S \setminus I$. Lemma \ref{lem 2} implies that $a_{1} \succ S \setminus I$. Consider $G + h'a_{1}$. By the same arguments as Case 1, we have a contradiction.
\vskip 5 pt

\indent Hence $\alpha_{1} \geq 2$ or $\alpha_{1} = 1$ but $\{a_{i}\} = I_{1}$ for some $i \in \{1, 2, ..., p - 1\}$. It is not difficult to see that there exists $j \in \{1, 2, ..., p - 2\}$ such that $a_{j} \in H_{l}$ and $a_{j + 1} \in H_{3 - l}$ for some $l \in \{1, 2\}$. Since $a_{j} \in H_{l}$, $x_{j} \succ H_{3 - l} \setminus \{a_{j + 1}\}$ by Lemma \ref{lem 2}. Similarly, $x_{j + 1} \succ H_{l} \setminus \{a_{j + 2}\}$ (it is possible that $a_{j + 2} \notin H_{l}$) because $a_{j + 1} \in H_{3 - l}$. Therefore, $\{x_{j}, x_{j + 1}\} \succ_{c} G - S$. Since $x_{j} \succ S \setminus I$, by Claim 1, $x_{j} \succ S$. Hence $\{x_{j}, x_{j + 1}\} \succ_{c} G$ contradicting $\gamma_{c}(G) = 3$. Thus Case 2 cannot occur and this completes the proof.
\qed

\indent We now have these following results as consequences.
\vskip 5 pt

\begin{thm}\label{thm MP1}
Let $G$ be a $3$-$\gamma_{c}$-critical graph. If $\kappa + 1 \leq \alpha$, then $G \in \mathcal{G}_{2}(3, 3)$ or $G$ satisfies $\kappa = \delta$.
\end{thm}
\proof
Suppose that $\kappa + 1 \leq \alpha$ and $G \notin \mathcal{G}_{2}(3, 3)$. By Theorem \ref{thm PM1}, $|S \setminus I| \geq 2$ or $\kappa \geq \delta$. Because $\kappa \leq \delta$ for any graph, it suffices to consider the case when $|S \setminus I| \geq 2$. Thus, Lemma \ref{lem PM2} yields that for a minimum cut set $S$, $G - S$ contains at least one singleton component $H$. Let $\{h\} = V(H)$. Thus $N_{G}(h) = S$. This implies that $\kappa = |S| = deg_{G}(h) \geq \delta$. Therefore $\kappa = \delta$ and this completes the proof.
\qed

\begin{thm}\label{thm MP2}
Let $G$ be a $3$-$\gamma_{c}$-critical graph such that $G \notin \mathcal{G}_{2}(3, 3)$. Then $\alpha = \kappa + q$ if and only if $\alpha = \delta + q$ for all $1 \leq q \leq 2$.
\end{thm}
\proof
We first suppose that $\alpha = \kappa + q$. By Theorem \ref{thm MP1}, $\kappa = \delta$ because $G \notin \mathcal{G}_{2}(3, 3)$. Thus $\alpha = \delta + q$.
\vskip 5 pt

\indent Conversely, let $\alpha = \delta + q$. Thus $\alpha = \delta + q \geq \kappa + q$. We first consider the case when $q = 2$. Therefore $\alpha \geq \kappa + 2$. By Theorem \ref{thm PM0}, $\kappa + 2 \geq \alpha$. Thus $\alpha = \kappa + 2$. We now consider the case when $q = 1$. Thus $\alpha = \delta + 1 \geq \kappa + 1$. If $\alpha = \kappa + 2$, then, by Theorem \ref{thm MP1}, $\kappa = \delta$. Thus $\alpha = \delta + 2$ contradicting $\alpha = \delta + 1$ and this completes the proof.
\qed

\indent We next show that there are at least two classes of $3$-$\gamma_{c}$-critical graphs with $\alpha \leq \kappa$ but $\kappa < \delta$. That is the condition $\kappa + 1 \leq \alpha$ in Theorems \ref{thm MP1} and \ref{thm MP2} is best possible. Let $M_{1}, M_{2}, ..., M_{s}$ be disjoint vertex sets where $|M_{i}| = m_{i} \geq 2$ for all $1 \leq i \leq s$ and $a_{1}, a_{2}, ..., a_{s}$ be isolated vertices. A graph $G$ in the class $\mathcal{G}_{3}(m_{0}, m_{1}, m_{2}, ..., m_{s})$ can be constructed from $M_{1}, M_{2}, ..., M_{s}$ and $a_{1}, a_{2}, ..., a_{s}$ by adding edges according to the join operation :
\begin{itemize}
  \item $a_{i} \vee ((\cup^{s}_{j = 0}M_{j}) \setminus M_{i})$ for $1 \leq i \leq s$,
  \item adding edges so that $M_{0}$ is a clique,
  \item adding edges so that $\cup^{s}_{i = 1}M_{i}$ is a clique.
\end{itemize}
\noindent A graph $G$ is illustrated in Figure 2.
\vskip 20 pt

\setlength{\unitlength}{0.7cm}
\begin{picture}(0, 7.5)
%Big G3

%vertices
\put(7.5, 3){\circle*{0.2}}
\put(9.5, 3){\circle*{0.2}}
\put(13.5, 3){\circle*{0.2}}

%oval
\put(10.5, 6){\oval(8, 1)}
\put(10.5, 0){\oval(2, 1)}

%a1 up
\put(7.5, 3){\line(-1, 3){0.85}}
\put(7.5, 3){\line(5, 2){6.2}}
%a1 down
\put(7.5, 3){\line(2, -3){2}}
\put(7.5, 3){\line(1, -1){2.5}}

%a2 down
\put(9.5, 3){\line(-1, 1){2.5}}
\put(9.35, 3){\line(2, 1){5}}

\put(9.5, 3){\line(0, -1){2.7}}
\put(9.5, 3){\line(1, -2){1.23}}

%as up
\put(13.5, 3){\line(1, 3){0.85}}
\put(13.5, 3){\line(-5, 2){6.2}}
%as down
\put(13.5, 3){\line(-2, -3){2}}
\put(13.5, 3){\line(-1, -1){2.5}}

%partition
\put(8.5, 5.5){\line(0, 1){1}}
\put(10.5, 5.5){\line(0, 1){1}}
\put(12.5, 5.5){\line(0, 1){1}}

\multiput(7.4, 3)(0,0.7){4}{\line(0,1){0.4}}
\multiput(7.6, 3)(0,0.7){4}{\line(0,1){0.4}}

\multiput(9.4, 3)(0,0.7){4}{\line(0,1){0.4}}
\multiput(9.6, 3)(0,0.7){4}{\line(0,1){0.4}}

\multiput(13.4, 3)(0,0.7){4}{\line(0,1){0.4}}
\multiput(13.6, 3)(0,0.7){4}{\line(0,1){0.4}}

%%%%%article%%%
\put(11.2, 3){$...$}
\put(11.2, 5.8){$...$}

\put(6.8, 3){$a_{1}$}
\put(8.8, 3){$a_{2}$}
\put(12.5, 3){$a_{s}$}

\put(7.5, 5.7){$M_{1}$}
\put(9.5, 5.7){$M_{2}$}
\put(13, 5.7){$M_{s}$}
\put(10.2, -0.2){$M_{0}$}
\put(5.4, -1.5){\footnotesize\textbf{Figure 2} : A graph $G$ in the class $\mathcal{G}_{3}(m_{0}, m_{1}, ..., m_{s})$.}
%\put(17, 6){\footnotesize $a_{0}$}
%%%%%%%END Big G3%%%%%%%
\end{picture}
\vskip 70 pt

\begin{lem}\label{lem g3}
If $G \in \mathcal{G}_{3}(m_{0}, m_{1}, ..., m_{s})$, then $G$ is a $3$-$\gamma_{c}$-critical graph satisfying $\alpha \leq \kappa < \delta$.
\end{lem}
\proof Let $m' \in M_{0}$. Therefore $\{m', a_{1}, a_{2}\} \succ_{c} G$. Thus $\gamma_{c}(G) \leq 3$. Suppose there exist $u, v \in V(G)$ such that $\{u, v\} \succ_{c} G$. Suppose that $\{u, v\} \cap \{a_{1}, a_{2}, ..., a_{s}\} \neq \emptyset$. By connected, $|\{u, v\} \cap \{a_{1}, a_{2}, ..., a_{s}\}| = 1$. Without loss of generality, let $u = a_{i}$ for some $i \in \{1, 2, ..., s\}$. To dominate $M_{i}$ and by connected, $v \in M_{j}$ where $j \in \{1, 2, ..., s\} \setminus \{i\}$. Thus $\{u, v\}$ does not dominate $a_{j}$. Therefore $\{u, v\} \cap \{a_{1}, a_{2}, ..., a_{s}\} = \emptyset$. By connected, either $\{u, v\} \subseteq M_{0}$ or $\{u, v\} \subseteq \cup^{s}_{i = 1}M_{i}$. Hence $\{u, v\}$ does not dominate $M_{1}$ if $\{u, v\} \subseteq M_{0}$ and $\{u, v\}$ does not dominate $M_{0}$ otherwise. Hence $\gamma_{c} \geq 3$ which gives $\gamma_{c}(G) = 3$.
\vskip 5 pt

\indent We next establish the criticality. For a pair of non-adjacent vertices $u, v \in V(G)$. Let $i \in \{1, 2, ..., s\}$ and $m_{i} \in M_{i}$. If $\{u, v\} = \{m_{0}, m_{i}\}$, then $D_{uv} = \{m_{0}, m_{i}\}$. If $\{u, v\} = \{a_{i}, m_{i}\}$ or $\{u, v\} = \{a_{i}, a_{j}\}$ for $j \neq i$, then $D_{uv} = \{a_{j}, m_{i}\}$. Thus $G$ is a $3$-$\gamma_{c}$-critical graph.
\vskip 5 pt

\indent We see that $G$ has $\{a_{1}, a_{2}, ..., a_{s}\}$ the minimum cut set as well as the maximum independent set. So $\alpha = \kappa = |S|$. It is easy to see that $\delta \geq s + 1$. Thus $\kappa < \delta$ and  this completes the proof.
\qed
\vskip 5 pt

\indent Consider a graph $G \in \mathcal{G}_{1}(b_{0}, .., b_{s - 1})$ when $b_{i} \geq 2$ for all $0 \leq i \leq s - 1$. A set $\{a_{0}, a_{1}, ...,$ $ a_{s - 1}, x\}$ is a minimum cut set of $G'$ with $s + 1$ vertices and a set $\{a_{0}, a_{1}, ..., a_{s}\}$ is the maximum independent set with $s + 1$ vertices. Hence, $\alpha = \kappa$. It is easy to see that $\delta \geq s + 1$. Therefore, $\delta > \kappa$.
%%%%%%%%%%%%%%%%%%%%%%%%%%%%%%%%%%%%%%%%%%%%%%%%%%%%%%%%%%%%%%%%%%%%%%%%%%%%%%%%%%%%%%%%%%%%%%%%%%%%%%%%%
\section{Discussion}
We conclude this paper by considering the Hamiltonian connected of $3$-$\gamma_{c}$-critical graphs according to $\alpha, \kappa$ and $\delta$. We also have the following result as a consequence of Theorem \ref{thm MP1}.

\begin{thm}\label{thm MP3}
Let $G$ be a $3$-$\gamma_{c}$-critical graph with $\delta \geq 2$. If $\kappa \leq \delta - 1$, then $\alpha \leq \kappa$ or $G \in \mathcal{G}_{2}(3, 3)$.
\end{thm}
\proof
Suppose that $\kappa \leq \delta - 1$. By Theorem \ref{thm MP1}, $G \in \mathcal{G}_{2}(3, 3)$ or $\alpha < \kappa + 1$. If $G \notin \mathcal{G}_{2}(3, 3)$, then $\alpha \leq \kappa + 1 - 1$. Therefore $\alpha \leq \kappa$ and this completes the proof.
\qed

\indent We now consider the graph $G \in \mathcal{G}_{1}(b_{0}, .., b_{s - 1})$ when $b_{i} = 1$ for all $0 \leq i \leq s - 1$ and $s \geq 4$. For example, the graph $G \in \mathcal{G}_{1}(1, 1, 1, 1)$ when $s = 4$ is illustrated in Figure 3 with $B_{0} = \{y\}$.
\vskip 10 pt

\setlength{\unitlength}{0.8cm}
\begin{center}
\begin{picture}(1, 8)

\put(0, 0){\circle*{0.2}}
\put(0, -0.5){\footnotesize{$a_{4}$}}

\put(-2, 2){\circle*{0.2}}

\put(0, 2){\circle*{0.2}}

\put(2, 2){\circle*{0.2}}

\put(0, 4){\circle*{0.2}}

\put(-2, 4){\circle*{0.2}}
\put(-2.5, 4.25){\footnotesize{$a_{0}$}}

\put(0, 7){\circle*{0.2}}
\put(0, 7.25){\footnotesize{$y$}}

\put(2, 7){\circle*{0.2}}
\put(2, 7.25){\footnotesize{$x$}}

\put(4, 4){\circle*{0.2}}
\put(4.25, 4.25){\footnotesize{$a_{3}$}}

\put(2, 4){\circle*{0.2}}

\put(0, 0){\line(0, 2){2}}
\put(0, 2){\line(-2, 0){2}}
\put(0, 0){\line(-1, 1){2}}
\put(-2, 2){\line(0, 1){2}}
\put(0, 2){\line(0, 1){2}}
\put(0, 2){\line(2, 5){2}}
\put(0, 0){\line(1, 1){2}}
\put(0, 2){\line(1, 0){2}}
\put(2, 2){\line(0, 1){2}}
\qbezier(-2, 2)(0, 1)(2, 2)
\put(0, 4){\line(0, 1){3}}
\put(2, 4){\line(0, 1){3}}
\put(-2, 2){\line(2, 1){4}}
\put(-2, 2){\line(3, 1){6}}
\qbezier(-2, 2)(-0.75, 4.75)(2, 7)
\put(0, 2){\line(-1, 1){2}}
\put(0, 2){\line(2, 1){4}}
\put(2, 2){\line(-2, 1){4}}
\put(2, 2){\line(-1, 1){2}}
\qbezier(2, 2)(1, 4.5)(2, 7)
\put(-2, 4){\line(2, 3){2}}
\put(-2, 4){\line(4, 3){4}}
\put(0, 4){\line(2, 3){2}}
\put(2, 4){\line(-2, 3){2}}
\put(4, 4){\line(-4, 3){4}}
\put(4, 4){\line(-2, 3){2}}
\put(-2, -1.5){\footnotesize\textbf{Figure 3} : $G \in \mathcal{G}_{1}(1, 1, 1, 1)$.}
%\put(17, 6){\footnotesize $a_{0}$}
\end{picture}
\end{center}
\vskip 40 pt

\noindent Clearly, $G$ satisfies $\kappa = \delta \geq 3$. It is well known that if $G$ is Hamiltonian connected, then $\omega(G - S) \leq |S| - 1$ for all cut set $S$ of $G$ where $\omega(G - S)$ is the number of components of $G - S$. We may assume that $B_{i} = \{x_{i}\}$ for all $1 \leq i \leq s - 1$ and $B_{0} = \{y\}$. If we let $S = \{y, x, x_{1}, x_{2}, ..., x_{s - 1}\}$, then $|S| = s + 1$ and $G - S$ has $a_{0}, a_{1}, ..., a_{s}$ as the $s + 1$ components. Thus, $\omega(G - S) = |S|$ and this implies that $G$ is not Hamiltonian connected. Hence, when $\kappa = \delta \geq 3$, a $3$-$\gamma_{c}$-critical graph does not need to be Hamiltonian connected. Thus, we may consider Hamiltonian connected property of these $3$-$\gamma_{c}$-critical graphs when $\kappa \leq \delta - 1$. By Theorem \ref{thm MP3}, $\alpha \leq \kappa$ because $\kappa \geq 3$. We note by Theorem \ref{thm CE} that every graph is Hamiltonian connected if $\alpha < \kappa$. Thus, the question that arises is :
\vskip 5 pt

\indent Is every $3$-$\gamma_{c}$-critical graph with $3 \leq \alpha = \kappa \leq \delta - 1$ Hamiltonian connected?
\vskip 15 pt

\noindent \textbf{Acknowledgement :} The authors express sincere thanks to Thiradet Jiarasuksakun for his valuable suggestion.

\end{document}